\documentclass[11pt]{article}

\usepackage[T1]{fontenc}
\usepackage[utf8]{inputenc}
\usepackage{lmodern}
\usepackage[margin=1in]{geometry}
\usepackage{microtype}

\usepackage{amsmath,amssymb,amsthm,mathtools}
\usepackage{bm}

\usepackage{graphicx}
\graphicspath{{figures/}}
\usepackage{booktabs}
\usepackage{array}

\usepackage[shortlabels]{enumitem}
\usepackage[numbers,sort&compress]{natbib}
\usepackage{xcolor}
\definecolor{referenceblue}{HTML}{194C70}
\usepackage[colorlinks=true,allcolors=referenceblue]{hyperref}
\usepackage[capitalise,noabbrev]{cleveref}

\theoremstyle{plain}
\newtheorem{theorem}{Theorem}[section]
\newtheorem{lemma}[theorem]{Lemma}
\newtheorem{proposition}[theorem]{Proposition}
\newtheorem{corollary}[theorem]{Corollary}

\theoremstyle{definition}
\newtheorem{definition}[theorem]{Definition}

\theoremstyle{remark}

\newcommand{\R}{\mathbb{R}}

\DeclareMathOperator{\diag}{diag}
\newcommand{\one}{\mathbf{1}}
\newcommand{\norm}[1]{\left\lVert#1\right\rVert}
\newcommand{\simplex}{\Delta_m}
\newtheorem{assumption}[theorem]{Assumption}
\newenvironment{procedure}[1]{\par\medskip\noindent\begin{minipage}{\linewidth}\textbf{#1}\par\smallskip}{\end{minipage}\par\medskip}
\allowdisplaybreaks[1]
\title{Nash Gradient Aggregation as Risk Budgeting}
\author{Jinze Zhao\\University of California, San Diego\\\texttt{jiz419@ucsd.edu}}
\date{}

\begin{document}

\maketitle

\begin{abstract}
Nash multi-task learning combines task gradients into a direction that improves every task. Risk budgeting allocates a portfolio so that each asset contributes a prescribed share of its risk. We show how to convert the weights of one problem into those of the other by using the matrix of gradient inner products as the portfolio covariance matrix. A unique positive solution exists exactly when the tasks admit a strict common descent direction, including cases with linearly dependent gradients.

How accurately must we solve for the weights before taking a training step? A cyclic coordinate solver updates one weight at a time and stops when every relative budget error is within a fixed tolerance smaller than one. We bound the number of sweeps needed for a fixed subproblem. For smooth, lower-bounded objectives, we prove that the resulting normalized update decreases every task objective for a specified step size and that the best squared Pareto stationarity measure over $T$ steps is $O(1/T)$. We also quantify the effects of using an older or regularized Gram matrix and of estimating the gradients themselves.

\end{abstract}

\newpage
\tableofcontents
\newpage

\section{Introduction}\label{sec:introduction}

How much weight should each task receive in a shared training step? A task with a large gradient can dominate the update. Two tasks with gradients that point in opposite directions can cancel each other. Nash multi-task learning (Nash-MTL) chooses the weights by bargaining over the improvement that each task receives. Portfolio risk budgeting asks a similar question about assets. A portfolio assigns each asset a weight so that its contribution to total risk matches a prescribed share.

To see the connection, let $g_i$ be the gradient of task $i$, let $G$ be the matrix with columns $g_i$, and let $a_i>0$ be the weight on that gradient. The proposed parameter update moves along $-v$, where $v=Ga$. Task $i$ then receives a first-order improvement proportional to $g_i^\top v$. Its weight times that improvement per unit step is
\[
 a_i g_i^\top v
 = a_i g_i^\top Ga
 = a_i(Qa)_i,\qquad Q=G^\top G .
\]
The last expression also appears in portfolio optimization. If $Q$ is an asset covariance matrix, then $a_i(Qa)_i$ is asset $i$'s contribution to the variance of the unnormalized holdings $a$. Thus the same product has an interpretation on each side. Prescribing its value $b_i>0$ for every $i$ gives
\begin{equation}\label{eq:central}
 Qa=b/a,\qquad a\in\R_{++}^m ,
\end{equation}
where division is componentwise. Equal budgets recover the original Nash-MTL equation and, on the portfolio side, equal risk contribution.

\paragraph{Computing a training step.}
Equation~\eqref{eq:central} is the first-order condition of a convex optimization problem from risk budgeting,
\begin{equation}\label{eq:potential-intro}
 \min_{a>0}\ \Phi_{Q,b}(a),\qquad
 \Phi_{Q,b}(a)=\frac12a^\top Qa-\sum_{i=1}^m b_i\log a_i .
\end{equation}
A portfolio optimizer can therefore compute the weights for a training step. Training also poses a question of its own. The gradients will change at the next iterate, so how accurately should we solve today's problem? We measure the relative error in each budget and stop when the largest error meets a tolerance. We prove that this test suffices for a step that decreases every smooth task objective, and we bound the number of coordinate sweeps needed to pass it.

\paragraph{Relation to existing work.}
\citet{navon2022} derive the equal-budget Nash-MTL equation and solve it by sequential convex optimization. For portfolio weights, \citet{spinu2013} gives the quadratic logarithmic objective in \eqref{eq:potential-intro}. \citet{griveau2013} develop a coordinate method for the volatility formulation, and \citet{palomarpackage} describe coordinate updates for the quadratic formulation used here. On the learning side, \citet{li2026toppo} write down the quadratic logarithmic program, and \citet{hu2026} analyze normalized Nash aggregation without requiring linearly independent gradients. We build on these results by tracing the error from an individual coordinate sweep through to the decrease in each task objective.

\begin{samepage}
\paragraph{Roadmap.}
Section~\ref{sec:methods} starts with the two decisions the algorithms make and derives the task utilities and portfolio risk contributions. Section~\ref{sec:correspondence} matches these quantities, proves when a positive solution exists, and gives the change of scale that converts one set of weights into the other. Section~\ref{sec:solver} uses this equivalence to compute the weights one coordinate at a time. Its stopping rule is the input to Section~\ref{sec:outer}, which turns a measured budget error into a bound on progress toward Pareto stationarity.

Section~\ref{sec:perturb} asks how that guarantee changes when the weights come from an earlier step, a regularized matrix, or noisy gradients. Section~\ref{sec:research} combines the inner and outer bounds to count the work of training and lays out experiments for the cases where that cost remains uncertain. Section~\ref{sec:conclusion} closes the paper. For a numerical illustration of the change of scale, Appendix~\ref{app:examples} works through two tasks by hand.
\end{samepage}

\section{The two optimization methods}\label{sec:methods}

At a single training step, Nash-MTL chooses how much each task gradient contributes to the update. A portfolio manager faces a related question when choosing how much each asset contributes to portfolio risk. Both decisions depend on pairwise interactions. For tasks we measure them with gradient inner products. For assets we measure them with covariances. We derive the weight equation in each setting and identify the output those weights produce.

\subsection{Notation and the multi-objective problem}
Let $f_1,\ldots,f_m:\R^p\to\R$ be differentiable task objectives. Collecting their gradients at $x$ gives the matrices
\begin{equation}\label{eq:gram}
\begin{aligned}
 g_i &= \nabla f_i(x)                    &&\text{(task gradient)},\\
 G   &= [g_1,\ldots,g_m]\in\R^{p\times m} &&\text{(one column per task)},\\
 Q   &= G^\top G,\qquad Q_{ij}=g_i^\top g_j &&\text{(Gram matrix)}.
\end{aligned}
\end{equation}
For any vector $a$, $a^\top Qa=(Ga)^\top(Ga)=\norm{Ga}^2\ge0$, so $Q$ is symmetric positive semidefinite, written $Q\succeq0$. We use the Euclidean norm for vectors and its induced norm for matrices.

We use $a>0$ for a vector with strictly positive coordinates. We use $w$ for weights whose sum is one. The admissible set of such nonnegative weights and the budget notation are
\[
 \simplex=\{w\in\R^m:w\ge0,\ \one^\top w=1\},\qquad
 b_i>0,\quad B=\sum_{i=1}^m b_i,\quad b_{\min}=\min_i b_i .
\]
The unweighted Nash-MTL choice is $b=\one$, for which $B=m$.

A parameter update takes the form $x^+=x-\eta d$. We call $d$ a \emph{strict common descent direction} when $g_i^\top d>0$ for all $i$. To see the sign convention, fix a task and write
\[
 \lim_{\eta\downarrow0}\frac{f_i(x-\eta d)-f_i(x)}{\eta}
 =-g_i^\top d.
\]
If each inner product is positive, each limit is negative. There are only finitely many tasks, so a sufficiently small positive step decreases all their objectives.

\begin{definition}[Pareto stationarity measure]\label{def:chi}
The distance from zero to the convex hull of the task gradients is
\begin{equation}\label{eq:chi}
 \chi(x)=\min_{w\in\simplex}\norm{G(x)w}.
\end{equation}
A point is first-order Pareto stationary if $\chi(x)=0$, that is, if some convex combination of its task gradients vanishes.
\end{definition}

\subsection{Nash-MTL and bargaining over a descent direction}
The first-order gain of task $i$ per unit step is $u_i=g_i^\top d$. A task that stays at $x$ receives zero gain. Nash bargaining chooses a direction with positive gain for every task and balances the gains by maximizing their weighted product. For radius $r>0$, taking the logarithm of the product gives the equivalent problem
\begin{equation}\label{eq:bargain}
\begin{aligned}
 \underset{d}{\operatorname{maximize}} \quad
 &U_b(d)=\sum_{i=1}^m b_i\log(g_i^\top d)\\
 \text{subject to}\quad &\norm d\le r,\qquad g_i^\top d>0\quad(i=1,\ldots,m).
\end{aligned}
\end{equation}
Indeed, $\log(\prod_i u_i^{b_i})=\sum_i b_i\log u_i$, and the logarithm preserves order. The marginal gain in the bargaining objective from increasing $u_i$ is $b_i/u_i$. With the same budget, a task gaining less from the current direction receives greater marginal weight.

For $b=\one$, \citet{navon2022} express the solution through the positive equation $Qa=\one/a$. For general budgets, the corresponding equation and direction are
\begin{equation}\label{eq:canonical}
\begin{aligned}
 Qa_\star &= b/a_\star,                    &&\text{(canonical weights)},\\
 v_\star &= Ga_\star,\qquad \norm{v_\star}^2=B, &&\text{(canonical direction)},\\
 d_r &= (r/\sqrt B)\,v_\star               &&\text{(bargaining direction of radius $r$)}.
\end{aligned}
\end{equation}
Division is componentwise. Section~\ref{sec:correspondence} proves when the weights exist and why their direction maximizes \eqref{eq:bargain}.

\begin{procedure}{Algorithm 1. Canonical Nash-MTL with positive budgets}
\begin{enumerate}[leftmargin=*]
\item Choose an initial point $x_0$, budgets $b>0$, and step sizes $\eta_t>0$.
\item At $x_t$, compute each task gradient, collect them in $G_t$, and form $Q_t=G_t^\top G_t$.
\item If $\chi(x_t)=0$, return $x_t$. Otherwise find $a_t>0$ solving $Q_ta_t=b/a_t$.
\item Form $v_t=G_ta_t$ and update $x_{t+1}=x_t-\eta_t v_t$. Return to Step 2.
\end{enumerate}
\end{procedure}

\subsection{The original Nash-MTL numerical procedure}
Step 3 asks for the weights at a fixed $Q$. For the original equal-budget case, \citet{navon2022} introduce $z_i(a)=(Qa)_i$ and work on the domain $a>0$, $Qa>0$. Taking logarithms on this positive domain gives
\[
 (Qa)_i=1/a_i
 \quad\Longleftrightarrow\quad
 a_i z_i(a)=1
 \quad\Longleftrightarrow\quad
 \underbrace{\log a_i+\log z_i(a)}_{\phi_i(a)}=0.
\]
Write $\phi(a)=\sum_i\phi_i(a)$. Each $\phi_i$ is concave because it is the sum of logarithms of positive affine functions. Its superlevel set $\{\phi_i\ge0\}$ is convex. The initial convex problem is
\begin{equation}\label{eq:ccp-initial}
\begin{aligned}
 &\min_{a\in\mathcal F}\;\one^\top Qa,\\
 &\mathcal F=\{a>0:Qa>0,\ \phi_i(a)\ge0\ \text{for all }i\}.
\end{aligned}
\end{equation}
Starting from a feasible $a^{(k)}$, the sequential convex procedure keeps this feasible set and uses the tangent of $\phi$ at the current iterate. The next subproblem is
\begin{equation}\label{eq:ccp-step}
 \min_{a\in\mathcal F}\;
 \one^\top Qa+\phi(a^{(k)})+
 \nabla\phi(a^{(k)})^\top(a-a^{(k)}).
\end{equation}
For coordinate $j$, the chain rule gives
\begin{align*}
 \partial_j\phi(a)
 &=\frac1{a_j}+\sum_i\frac{\partial_j(Qa)_i}{(Qa)_i}\\
 &=\frac1{a_j}+\sum_i\frac{Q_{ij}}{(Qa)_i}.
\end{align*}
In vector form, $\nabla\phi(a)=a^{-1}+Q^\top(Qa)^{-1}$. The risk-budgeting formulation below gives a second way to obtain the weights for Step 3.

\subsection{Risk budgeting and the contribution of each asset}
Let $\Sigma\succeq0$ be a covariance matrix and $w\in\simplex$ a portfolio. Its variance is $V(w)=w^\top\Sigma w$, and its volatility is $\sigma(w)=\sqrt{V(w)}$. For $V(w)>0$, the volatility change from increasing one asset weight is
\begin{align*}
 \frac{\partial\sigma(w)}{\partial w_i}
 &=\frac{1}{2\sqrt{V(w)}}\frac{\partial V(w)}{\partial w_i}
       &&\text{(chain rule)}\\
 &=\frac{2(\Sigma w)_i}{2\sqrt{V(w)}}
       &&\text{($\Sigma$ is symmetric)}\\
 &=\frac{(\Sigma w)_i}{\sigma(w)}.
\end{align*}
Multiplying this marginal change by the amount held defines the asset's risk contribution, $\operatorname{RC}_i(w)=w_i(\Sigma w)_i/\sigma(w)$. Summing over assets recovers the full portfolio volatility, since
\begin{align*}
 \sum_i\operatorname{RC}_i(w)
 &=\frac{\sum_i w_i(\Sigma w)_i}{\sigma(w)}\\
 &=\frac{w^\top\Sigma w}{\sigma(w)}
  =\sigma(w).
\end{align*}
Risk budgeting prescribes the fraction $\beta_i=b_i/B$ of that total attributable to each asset. Multiplying the target equality by $\sigma(w)>0$ gives
\begin{equation}\label{eq:rb}
\begin{aligned}
 \operatorname{RC}_i(w)=\beta_i\sigma(w)
 &\quad\Longleftrightarrow\quad w_i(\Sigma w)_i=\beta_i\sigma(w)^2\\
 &\quad\Longleftrightarrow\quad w_i(\Sigma w)_i=\frac{b_i}{B}w^\top\Sigma w.
\end{aligned}
\end{equation}
Equal budgets give each asset the same share of volatility. This is risk parity.

\citet{spinu2013} solves for positive auxiliary weights first, then normalizes them to obtain the portfolio according to
\begin{equation}\label{eq:spinu}
 \min_{a>0}\ \Phi_{\Sigma,b}(a)
 =\frac12a^\top\Sigma a-\sum_i b_i\log a_i,
 \qquad w=\frac{a}{\one^\top a}.
\end{equation}
Differentiation gives
\begin{equation}\label{eq:derivatives}
\begin{aligned}
 \partial_i\Phi_{\Sigma,b}(a)&=(\Sigma a)_i-b_i/a_i,\\
 \nabla\Phi_{\Sigma,b}(a)&=\Sigma a-b/a,\\
 \nabla^2\Phi_{\Sigma,b}(a)&=\Sigma+\diag(b_i/a_i^2).
\end{aligned}
\end{equation}
At an interior minimizer, the gradient vanishes. The resulting equation $\Sigma a=b/a$ is precisely the Nash weight equation \eqref{eq:canonical} with $\Sigma$ in place of $G^\top G$.

\begin{procedure}{Algorithm 2. Quadratic risk budgeting by damped Newton steps}
\begin{enumerate}[leftmargin=*]
\item Supply $\Sigma\succeq0$ with $\ker(\Sigma)\cap\R_+^m=\{0\}$, budgets $b>0$, a positive initial vector $a$, and a residual tolerance $\delta>0$.
\item Compute $q=\Sigma a$. If $\max_i|a_iq_i/b_i-1|\le\delta$, return the portfolio $w=a/(\one^\top a)$.
\item Find the Newton direction $p$ by solving
\[
 [\Sigma+\diag(b_i/a_i^2)]p=-(q-b/a).
\]
\item Choose a positive step $\tau$ by the backtracking rule in Appendix~\ref{app:newton}. Set $a\leftarrow a+\tau p$ and return to Step 2.
\end{enumerate}
\end{procedure}

The coordinate method in Section~\ref{sec:solver} minimizes the same quadratic objective. The risk-budgeting literature also uses the volatility objective $\sqrt{a^\top\Sigma a}-\sum_i b_i\log a_i$, whose coordinate derivative includes the volatility in its denominator \citep{griveau2013}. The quadratic coordinate update used here appears explicitly in \citet{palomarpackage}.

\section{Exact correspondence and its domain of validity}\label{sec:correspondence}

When can the weights computed for a portfolio also give a bargaining solution? A positive solution to the weight equation exists exactly when the tasks admit a strict common descent direction. Multiplication by $G$ then gives the Nash direction. Dividing the same weights by their sum gives the portfolio. We prove both constructions and show how the portfolio's quadratic risk recovers the original scale.

\subsection{Common descent, the nonnegative kernel, and existence}
\begin{lemma}[Three descriptions of common descent]\label{lem:alternative}
For a fixed gradient matrix $G$, the following are equivalent.
\begin{enumerate}[label=(\roman*),leftmargin=*]
\item The gradient hull has positive distance from zero, $\chi>0$.
\item The only nonnegative vector in the nullspace of $G$ is zero.
\item There is an $h$ such that $g_i^\top h>0$ for all $i$.
\end{enumerate}
\end{lemma}
\begin{proof}
\emph{Equivalence of (i) and (ii).}
The simplex is compact and $w\mapsto\norm{Gw}$ is continuous, so the minimum defining $\chi$ is attained. If $\chi=0$, any minimizing $w$ satisfies
\[
 w\ge0,\qquad \one^\top w=1,\qquad Gw=0.
\]
This $w$ is a nonzero nonnegative null vector. Conversely, given any $z\ge0$, $z\ne0$, with $Gz=0$, set $s=\one^\top z>0$. Its normalization has the properties
\[
 w=z/s\ge0,\qquad \one^\top w=1,\qquad Gw=(Gz)/s=0.
\]
Thus $\chi=0$ exactly when such a $z$ exists.

\emph{From (i) to (iii).}
Let $\bar w$ minimize $\norm{Gw}^2/2$ over the simplex and set $h=G\bar w$. Move from $\bar w$ toward the $i$th simplex vertex $e_i$. For $0\le\tau\le1$, the resulting gradient mixture is
\[
 G[(1-\tau)\bar w+\tau e_i]=h+\tau(g_i-h).
\]
Along this feasible segment the squared norm is smallest at $\tau=0$. Its right derivative there satisfies
\begin{align*}
 0
 &\le \left.\frac{d}{d\tau}\frac12\norm{h+\tau(g_i-h)}^2\right|_{\tau=0^+}\\
 &=h^\top(g_i-h)\\
 &=g_i^\top h-\norm h^2.
\end{align*}
Hence $g_i^\top h\ge\norm h^2=\chi^2>0$ for every $i$.

\emph{From (iii) to (i).}
Let $\alpha=\min_i g_i^\top h>0$. Every $w\in\simplex$ satisfies
\begin{align*}
 \norm h\,\norm{Gw}
 &\ge h^\top Gw &&\text{(Cauchy--Schwarz)}\\
 &=\sum_i w_i g_i^\top h\\
 &\ge\alpha\sum_i w_i=\alpha .
\end{align*}
Since $h\ne0$, minimizing over $w$ gives $\chi\ge\alpha/\norm h>0$.
\end{proof}

\begin{theorem}[Existence and uniqueness, including singular Gram matrices]\label{thm:exist}
Let $Q=G^\top G$ and $b>0$.
\begin{enumerate}[label=(\alph*),leftmargin=*]
\item If $\chi>0$, every nonempty finite sublevel set of $\Phi_{Q,b}$ is compact in $\R_{++}^m$. The function has a unique minimizer $a_\star>0$, characterized by $Qa_\star=b/a_\star$.
\item If $\chi=0$, then $\inf_{a>0}\Phi_{Q,b}(a)=-\infty$ and the weight equation has no positive solution.
\end{enumerate}
\end{theorem}
\begin{proof}
\emph{Bound the weights from above.}
Assume $\chi>0$. Write $s=\one^\top a$ and $w=a/s$. Then $w$ belongs to the simplex, so
\begin{align*}
 a^\top Qa &=\norm{Ga}^2\\
           &=s^2\norm{Gw}^2\\
           &\ge s^2\chi^2 .
\end{align*}
Since each $a_i\le s$ and each $b_i>0$, we also have $\sum_i b_i\log a_i\le B\log s$. Combining the two bounds gives
\begin{equation}\label{eq:coercive}
 \Phi_{Q,b}(a)\ge \frac{\chi^2}{2}s^2-B\log s=:\psi(s).
\end{equation}
Fix a finite level $C$ with a nonempty sublevel set $\mathcal S_C$. The function $\psi(s)$ tends to $+\infty$ as $s\to\infty$, since $\chi>0$. Choose $R>0$ large enough that $\psi(s)>C$ for every $s>R$. Then
\[
 a\in\mathcal S_C
 \quad\Longrightarrow\quad
 \psi(s)\le\Phi_{Q,b}(a)\le C
 \quad\Longrightarrow\quad
 a_i\le s\le R\quad\text{for all }i .
\]

\emph{Bound each weight away from zero.}
Isolating the logarithm of $a_i$ for $a\in\mathcal S_C$ gives
\begin{align*}
 C
 &\ge \Phi_{Q,b}(a)\\
 &\ge -b_i\log a_i-\sum_{j\ne i}b_j\log a_j
      &&(a^\top Qa\ge0)\\
 &\ge -b_i\log a_i-(B-b_i)\log R
      &&(a_j\le R).
\end{align*}
Rearrange, divide by $b_i>0$, and exponentiate to get
\begin{equation}\label{eq:lower-coord}
 a_i\ge\ell_i:=\exp\!\left(-\frac{C+(B-b_i)\log R}{b_i}\right)>0.
\end{equation}
Thus $\mathcal S_C\subseteq\prod_i[\ell_i,R]$. The box is compact and lies in the positive orthant. Continuity of $\Phi_{Q,b}$ makes its sublevel set closed in the box, hence compact.

\emph{Attain the minimum and identify it.}
Choose any $a^{(0)}>0$ and take $C=\Phi_{Q,b}(a^{(0)})$. The continuous objective attains its minimum at some $a_\star\in\mathcal S_C$. For every point outside this set, $\Phi_{Q,b}(a)>C\ge\Phi_{Q,b}(a_\star)$, so $a_\star$ is a global minimizer. For each nonzero $u$,
\begin{align*}
 u^\top\nabla^2\Phi_{Q,b}(a)u
 &=u^\top Qu+\sum_i\frac{b_i}{a_i^2}u_i^2\\
 &=\norm{Gu}^2+\sum_i\frac{b_i}{a_i^2}u_i^2\\
 &>0 .
\end{align*}
At least one $u_i$ is nonzero and its coefficient is positive. Consequently the objective has positive second derivative along every nonconstant segment and is strictly convex. Its minimizer is unique. On the open domain $\R_{++}^m$, the first-order condition reads
\[
 \nabla\Phi_{Q,b}(a_\star)=0
 \quad\Longleftrightarrow\quad
 Qa_\star-b/a_\star=0 .
\]
Conversely, for any zero $a$ of this gradient and any $y>0$, convexity gives $\Phi_{Q,b}(y)\ge\Phi_{Q,b}(a)+\nabla\Phi_{Q,b}(a)^\top(y-a)=\Phi_{Q,b}(a)$. Thus the weight equation characterizes the minimizer.

\emph{At $\chi=0$, follow a null vector.}
By Lemma~\ref{lem:alternative}, there is a $z\ge0$, $z\ne0$, with $Gz=0$. Along $a(t)=a^{(0)}+tz$, $t\ge0$, we have
\begin{align*}
 Ga(t)&=Ga^{(0)}+tGz=Ga^{(0)},\\
 \Phi_{Q,b}(a(t))
 &=\tfrac12\norm{Ga^{(0)}}^2-\sum_i b_i\log(a_i^{(0)}+tz_i)
 \longrightarrow-\infty .
\end{align*}
Every logarithm is nondecreasing in $t$, and at least one tends to $+\infty$. A positive solution of $Qa=b/a$ would also imply
\[
 0=(Gz)^\top Ga=z^\top Qa=\sum_i z_i b_i/a_i>0,
\]
where the strict inequality follows from $z\ge0$, $z\ne0$, and $b_i/a_i>0$. This contradiction excludes a positive solution.
\end{proof}

\subsection{From positive weights to the bargaining optimum}
Existence gives the weights. To show what the tasks gain from them, form the direction $v=Ga_\star$. Its inner product with each gradient is the corresponding entry of $Qa_\star$.

\begin{theorem}[Bargaining direction]\label{thm:bargain}
Assume $\chi>0$. Let $a_\star$ minimize $\Phi_{Q,b}$ and put $v=Ga_\star$. Then
\begin{equation}\label{eq:utilities}
 g_i^\top v=\frac{b_i}{a_{\star,i}}>0,\qquad
 \norm v^2=B.
\end{equation}
For every $r>0$, $d_r=(r/\sqrt B)v$ uniquely maximizes \eqref{eq:bargain}.
\end{theorem}
\begin{proof}
\emph{Compute the task gains and the direction length.} For each $i$,
\[
 g_i^\top v=g_i^\top Ga_\star=(Qa_\star)_i=b_i/a_{\star,i}.
\]
The sum of the weighted gains is
\begin{align*}
 \norm v^2 &= a_\star^\top Qa_\star\\
 &=\sum_i a_{\star,i}(Qa_\star)_i\\
 &=\sum_i a_{\star,i}\frac{b_i}{a_{\star,i}}=B .
\end{align*}
Set $c=r/\sqrt B$ and $d_r=cv$. Then $g_i^\top d_r=c\,b_i/a_{\star,i}>0$ and $\norm{d_r}=c\sqrt B=r$. This proves feasibility.

\emph{Compare with any other feasible direction $d$.}
The tangent inequality for the logarithm at $u_0>0$ is
$\log u-\log u_0\le(u-u_0)/u_0$ for $u>0$. Applying it task by task at $u_0=g_i^\top d_r$ with $u=g_i^\top d$ and summing with weights $b_i$ yields
\begin{align*}
 U_b(d)-U_b(d_r)
 &\le\sum_i\frac{b_i}{g_i^\top d_r}\,g_i^\top(d-d_r)\\
 &=\frac1c\sum_i a_{\star,i}g_i^\top(d-d_r)
      &&(g_i^\top d_r=c\,b_i/a_{\star,i})\\
 &=\frac1c v^\top(d-cv)
      &&\bigl(v=\textstyle\sum_i a_{\star,i}g_i\bigr)\\
 &\le\frac1c\norm v\,\norm d-\norm v^2
      &&\text{(Cauchy--Schwarz)}\\
 &\le\frac1c\sqrt B\,r-B=0
      &&(\norm d\le r).
\end{align*}
Any maximizer must attain equality in both of the last two inequalities. Thus $\norm d=r$ and $d$ is a positive multiple of $v$. These conditions fix $d=r v/\norm v=cv=d_r$.
\end{proof}

\subsection{The exact risk-budgeting map}
The bargain fixes the length of $Ga_\star$. A fully invested portfolio instead fixes the sum of its weights. Normalization moves between these conventions because both the individual contributions $a_i(Qa)_i$ and the total $a^\top Qa$ scale quadratically.

\begin{theorem}[Two-way correspondence]\label{thm:map}
Let $Q\succeq0$ with $\ker(Q)\cap\R_+^m=\{0\}$, and let $b>0$. The following three objects exist and are unique. Each determines the other two.
\begin{enumerate}[label=(\roman*),leftmargin=*]
\item the minimizer $a_\star$ of $\Phi_{Q,b}$,
\item the positive solution of $Qa=b/a$,
\item the positive simplex vector $w_\star$ satisfying
\begin{equation}\label{eq:portfolio-q}
 w_i(Qw)_i=(b_i/B)w^\top Qw\qquad(i=1,\ldots,m).
\end{equation}
\end{enumerate}
The conversion in each direction is
\begin{equation}\label{eq:inverse-map}
 w_\star=\frac{a_\star}{\one^\top a_\star},
 \qquad
 a_\star=\sqrt{\frac{B}{w_\star^\top Qw_\star}}\,w_\star .
\end{equation}
With $Q=G^\top G$, $Ga_\star$ is the canonical Nash direction. With $Q=\Sigma$, $w_\star$ is the risk-budgeting portfolio.
\end{theorem}
\begin{proof}
\emph{From the minimizer to the positive equation.}
Factor $Q=G^\top G$ using its positive semidefinite square root. Since $z^\top Qz=\norm{Gz}^2$, $Q$ and $G$ have the same kernel. The hypothesis and Lemma~\ref{lem:alternative} imply $\chi>0$. Theorem~\ref{thm:exist} now identifies (i) and (ii).

\emph{From canonical weights to portfolio weights.}
Set $a=a_\star$, $s=\one^\top a$, and $w=a/s$. Both the individual contributions and their sum acquire a factor $1/s^2$, since
\begin{align*}
 w_i(Qw)_i &=\frac{a_i(Qa)_i}{s^2}
            =\frac{b_i}{s^2},\\
 w^\top Qw &=\sum_i w_i(Qw)_i
            =\frac{B}{s^2}.
\end{align*}
The denominator is positive. Dividing the first line by the second shows that the $i$th contribution has share $(b_i/s^2)/(B/s^2)=b_i/B$, which is (iii).

\emph{From portfolio weights to canonical weights.}
Take any $w>0$ in the simplex satisfying \eqref{eq:portfolio-q}. Its quadratic risk $q=w^\top Qw=\norm{Gw}^2$ is positive by the kernel condition. Choose $c=\sqrt{B/q}$ and put $a=cw$. For each $i$,
\begin{align*}
 a_i(Qa)_i
 &=c^2 w_i(Qw)_i &&(a=cw)\\
 &=c^2(b_i/B)q &&\text{by \eqref{eq:portfolio-q}}\\
 &=(B/q)(b_i/B)q=b_i .
\end{align*}
Division by $a_i>0$ recovers (ii). It remains to check that we recover the original vector in either order. Starting with $w$, the normalized vector $cw/(\one^\top cw)$ equals $w$ since $\one^\top w=1$. Starting with $a$, the forward calculation gave $q=B/s^2$, hence $c=\sqrt{B/q}=s$ and $cw=s(a/s)=a$. The two maps are inverse, so uniqueness of the canonical solution gives uniqueness of the portfolio.
\end{proof}

\begin{center}
\begin{minipage}{\linewidth}
\centering\small
\refstepcounter{table}\textbf{Table~\thetable.} What the same equation represents in the two settings.\par\smallskip
\begin{tabular}{p{0.23\linewidth}p{0.32\linewidth}p{0.32\linewidth}}
\toprule
Quantity & Gradient aggregation & Risk budgeting\\\midrule
Interaction matrix $Q$ & Gradient Gram matrix $G^\top G$ & Asset covariance $\Sigma$\\
Budget $b_i$ & Task bargaining weight & Target risk share $b_i/B$\\
Positive weights $a$ & Solve $Qa=b/a$ & Minimize $\Phi_{Q,b}$\\
Simplex weights $w$ & Coefficients in $Gw$ & Asset holdings\\
$a_i(Qa)_i$ & Weighted task gain & Scaled contribution to variance\\
Output & Direction $Ga$ or $Gw$ & Portfolio $w$\\\bottomrule
\end{tabular}
\end{minipage}
\end{center}

\subsection{Scale invariance and the meaning of normalization}
\begin{proposition}[Effect of changing units or budgets]\label{prop:scale}
Assume $\chi>0$. Write $v=Ga_\star$ and $w_\star=a_\star/(\one^\top a_\star)$.
\begin{enumerate}[label=(\alph*),leftmargin=*]
\item Scaling each task gradient by $c_i>0$, so $G'=GD$ with $D=\diag(c_i)$, gives $a'_\star=D^{-1}a_\star$ and $G'a'_\star=v$.
\item Scaling every budget by the same $\gamma>0$ gives $a'_\star=\sqrt\gamma\,a_\star$, $w'_\star=w_\star$, and $Ga'_\star=\sqrt\gamma\,v$.
\end{enumerate}
\end{proposition}
\begin{proof}
\emph{Change of task units.} The new Gram matrix is $Q'=(GD)^\top GD=DQD$. For $a'=D^{-1}a_\star$ we obtain
\begin{align*}
 Q'a' &=DQD(D^{-1}a_\star)\\
      &=DQa_\star\\
      &=D(b/a_\star)=b/a'.
\end{align*}
The last equality holds coordinate by coordinate since $c_i b_i/a_{\star,i}=b_i/(a_{\star,i}/c_i)$. The matrix $D$ maps the nonnegative cone bijectively to itself, so $GD$ retains the nonnegative-kernel condition. Uniqueness identifies $a'$ as the new canonical solution, and $G'a'=GD(D^{-1}a_\star)=v$.

\emph{Change of budget scale.} For the scaled budgets $\gamma b$, the candidate $a'=\sqrt\gamma\,a_\star$ satisfies
\begin{align*}
 Qa' &=\sqrt\gamma\,Qa_\star\\
     &=\sqrt\gamma\,(b/a_\star)\\
     &=(\gamma b)/a'.
\end{align*}
Again uniqueness identifies the new weights. Normalizing and multiplying by $G$ give, respectively,
\[
 \frac{a'}{\one^\top a'}
 =\frac{\sqrt\gamma\,a_\star}{\sqrt\gamma\,\one^\top a_\star}
 =w_\star,\qquad
 Ga'=\sqrt\gamma\,v .
\]
\end{proof}

For a fixed $G$, put $s=\one^\top a_\star$. Then $Gw_\star=v/s$. Using $v\ne0$ to compare step sizes gives
\begin{align*}
 x-\widetilde\eta Gw_\star=x-\eta v
 &\quad\Longleftrightarrow\quad (\widetilde\eta/s)v=\eta v\\
 &\quad\Longleftrightarrow\quad \widetilde\eta=\eta s .
\end{align*}
Thus the two directions generate the same parameter step when their step sizes are related by $s$. For a fixed step size, changing $s$ changes the step length. This distinction matters as $G$ varies from one training step to the next. Section~\ref{sec:outer} studies the update using $Gw$.

\section{Transferring a coordinate solver}\label{sec:solver}
At a fixed parameter vector, Section~\ref{sec:correspondence} identifies the Nash weights with the minimizer of a risk-budgeting objective. How much work does it take to compute those weights? Holding all but one weight fixed reduces the problem to a quadratic equation with one positive root. Cycling through these scalar updates gives a solver. We derive its stopping rule and prove how many sweeps suffice to meet a prescribed tolerance.

\subsection{One coordinate at a time}
Fix $Q=G^\top G$ with $\chi>0$, and write $\Phi=\Phi_{Q,b}$. Each simplex vertex $e_i$ is feasible in the definition of $\chi$, so
\[
 Q_{ii}=\norm{Ge_i}^2\ge\min_{w\in\simplex}\norm{Gw}^2=\chi^2>0.
\]
To update the $i$th weight, hold the other coordinates fixed and set $u=a_i$. Write $r_i=\sum_{j\ne i}Q_{ij}a_j$. Grouping the terms containing $u$ and using $Q_{ji}=Q_{ij}$ gives
\begin{align*}
 \tfrac12a^\top Qa
 &=\tfrac12Q_{ii}u^2+\tfrac12u\sum_{j\ne i}(Q_{ij}+Q_{ji})a_j
    +\tfrac12\sum_{j,k\ne i}a_jQ_{jk}a_k\\
 &=\tfrac12Q_{ii}u^2+r_i u+\tfrac12\sum_{j,k\ne i}a_jQ_{jk}a_k.
\end{align*}
Dropping terms independent of $u$ leaves a one-variable minimization on $u>0$. Its first two derivatives are
\begin{align*}
 h_i(u)&=\tfrac12Q_{ii}u^2+r_i u-b_i\log u,\\
 h_i'(u)&=Q_{ii}u+r_i-b_i/u,\\
 h_i''(u)&=Q_{ii}+b_i/u^2>0.
\end{align*}
The derivative increases from $-\infty$ to $+\infty$, so its unique zero minimizes $h_i$. Multiplying the first-order equation by $u$ gives
\[
 h_i'(u)=0
 \quad\Longleftrightarrow\quad
 Q_{ii}u^2+r_i u-b_i=0.
\]
Write $D_i=\sqrt{r_i^2+4Q_{ii}b_i}$. Since $D_i>|r_i|$ and $Q_{ii}>0$, the negative-sign root is negative. The positive root is
\begin{equation}\label{eq:ccd-root}
 u=\frac{-r_i+D_i}{2Q_{ii}}.
\end{equation}
When $r_i\ge0$, subtracting $r_i$ from a nearby $D_i$ loses precision. Multiplying numerator and denominator by $D_i+r_i$ gives
\begin{align}
 u&=\frac{(D_i-r_i)(D_i+r_i)}{2Q_{ii}(D_i+r_i)}\nonumber\\
 &=\frac{4Q_{ii}b_i}{2Q_{ii}(D_i+r_i)}
   =\frac{2b_i}{D_i+r_i}.\label{eq:ccd-stable}
\end{align}
When $r_i<0$, the numerator in \eqref{eq:ccd-root} already adds two positive terms.

Store $q=Qa$ to reuse it across coordinates. Before the update, $r_i=q_i-Q_{ii}a_i$. Writing $\Delta_i=u-a_i$ for the change, the new matrix-vector product is
\begin{align}
 Q(a+\Delta_i e_i)&=Qa+\Delta_iQe_i\nonumber\\
 &=q+\Delta_i Q_{:i}.\label{eq:cache}
\end{align}
Only one column of $Q$ enters this calculation. Each update takes $O(m)$ operations, and one pass through all $m$ coordinates takes $O(m^2)$.

\subsection{A stopping rule measured in relative budget error}
At a solution, coordinate $i$ contributes exactly $b_i$. Measuring a trial vector's contribution in units of $b_i$ gives a common scale for all coordinates.
\begin{definition}\label{def:residual}
For $a>0$, the signed relative error and the largest absolute relative error are
\begin{align}
 e_i(a;Q,b)&=\frac{a_i(Qa)_i-b_i}{b_i},\nonumber\\
 \mathcal R(a;Q,b)&=\max_{1\le i\le m}|e_i(a;Q,b)|.\label{eq:residual}
\end{align}
\end{definition}
For instance, $\mathcal R\le0.1$ places each contribution between $0.9b_i$ and $1.1b_i$. The same condition measures accuracy in the Nash weight equation. It can be checked with the currently stored $a$ and $q$.

\begin{procedure}{Algorithm 3. Cyclic coordinate descent (CCD) for quadratic risk budgeting}
\begin{enumerate}[leftmargin=*]
\item Start with $Q\succeq0$ satisfying $\chi>0$, budgets $b>0$, weights $a>0$ and a tolerance $0<\delta<1$. Compute $q=Qa$.
\item If $\max_i|a_iq_i/b_i-1|\le\delta$, return $a$.
\item For $i=1,\ldots,m$, compute $r_i=q_i-Q_{ii}a_i$. Use \eqref{eq:ccd-stable} when $r_i\ge0$ and \eqref{eq:ccd-root} otherwise. With the resulting $u$, set
\[
 \Delta_i=u-a_i,\qquad a_i\leftarrow u,\qquad q\leftarrow q+\Delta_iQ_{:i}.
\]
\item At the end of the sweep, return to Step 2.
\end{enumerate}
\end{procedure}

One starting choice is $a_i=\sqrt{b_i/Q_{ii}}$, the solution for the diagonal part of $Q$. Successive outer iterations may also reuse the preceding positive weights. In floating-point arithmetic, periodically recompute $Qa$ and recompute it before accepting the final residual. Algorithm 3 uses the quadratic Spinu objective. The volatility objective in \citet{griveau2013} leads to a different coordinate equation.

\subsection{How much one sweep reduces the objective}
The objective decrease records how far the coordinates moved during a sweep. Those same movements determine the gradient left over at the end of the sweep. Bringing the two calculations together yields a contraction.

\begin{theorem}[Linear contraction per sweep]\label{thm:ccd}
Let $a^{(0)}>0$ be the starting point. Write $a_\star$ for the minimizer, $a^{(k)}$ for the vector after $k$ sweeps, and define
\[
 \mathcal S=\{a>0:\Phi(a)\le\Phi(a^{(0)})\},\qquad
 H(a)=\Phi(a)-\Phi(a_\star),\qquad H_k=H(a^{(k)}).
\]
Choose $R$ bounding every coordinate on $\mathcal S$. Let $U$ be the strictly upper triangular part of $Q$ in the order used for the sweep. Define
\[
 \mu=b_{\min}/R^2,\qquad
 \mu_c=\min_i(Q_{ii}+b_i/R^2),\qquad M=\norm U.
\]
If $M>0$, every complete sweep in exact arithmetic satisfies
\begin{equation}\label{eq:contraction}
 H_{k+1}\le\rho H_k,\qquad
 \rho=\frac{M^2}{M^2+\mu\mu_c}<1.
\end{equation}
Thus $H_k\le\rho^k H_0$ and $\norm{a^{(k)}-a_\star}^2\le 2H_k/\mu$. If $M=0$, one sweep suffices.
\end{theorem}

\begin{proof}
\emph{Curvature on the visited set.}
Coordinate minimization decreases or preserves $\Phi$, so all intermediate iterates belong to $\mathcal S$. This set is compact by Theorem~\ref{thm:exist} and convex by convexity of $\Phi$. For any $a\in\mathcal S$ and vector $v$,
\begin{align*}
 v^\top\nabla^2\Phi(a)v
 &=v^\top Qv+\sum_i b_i v_i^2/a_i^2\\
 &\ge\sum_i b_i v_i^2/R^2 && (Q\succeq0,\ a_i\le R)\\
 &\ge\mu\norm v^2. &&
\end{align*}
Take $a,z\in\mathcal S$ and set $v=z-a$. Integrating along the segment inside $\mathcal S$ gives
\begin{align*}
 \Phi(z)-\Phi(a)-\nabla\Phi(a)^\top v
 &=\int_0^1(1-t)v^\top\nabla^2\Phi(a+tv)v\,dt\\
 &\ge\mu\norm v^2\int_0^1(1-t)\,dt\\
 &=\tfrac\mu2\norm v^2.
\end{align*}
Consequently,
\begin{equation}\label{eq:strongconv}
 \Phi(z)\ge\Phi(a)+\nabla\Phi(a)^\top(z-a)+\tfrac\mu2\norm{z-a}^2.
\end{equation}
Setting $a=a_\star$ gives $H(z)\ge(\mu/2)\norm{z-a_\star}^2$, since $\nabla\Phi(a_\star)=0$. Setting $z=a_\star$ gives a bound in the other direction. Put $p=\nabla\Phi(a)$ and $h=a-a_\star$. Completing the square,
\begin{align}
 H(a)&\le p^\top h-\tfrac\mu2\norm h^2\nonumber\\
 &=\frac{\norm p^2}{2\mu}-\frac\mu2\norm{h-p/\mu}^2\nonumber\\
 &\le\frac{\norm{\nabla\Phi(a)}^2}{2\mu}.\label{eq:gap-gradient}
\end{align}

\emph{Decrease during the sweep.}
Let $z^{(i)}$ denote the vector just after updating coordinate $i$ and $z^{(0)}$ the starting vector. Write $s_i=z_i^{(i)}-z_i^{(i-1)}$. Along that coordinate segment, $h_i''(u)\ge Q_{ii}+b_i/R^2\ge\mu_c$. Apply \eqref{eq:strongconv} in one dimension, expanding at the new minimizer. Since $h_i'(z_i^{(i)})=0$,
\begin{align*}
 \Phi(z^{(i-1)})-\Phi(z^{(i)})
 &= h_i(z_i^{(i-1)})-h_i(z_i^{(i)})\\
 &\ge h_i'(z_i^{(i)})(-s_i)+\tfrac{\mu_c}{2}s_i^2\\
 &=\tfrac{\mu_c}{2}s_i^2.
\end{align*}
The intermediate objective values cancel in the sum, giving
\begin{align}
 D:=\Phi(z^{(0)})-\Phi(z^{(m)})
 &=\sum_{i=1}^m\bigl[\Phi(z^{(i-1)})-\Phi(z^{(i)})\bigr]\nonumber\\
 &\ge\tfrac{\mu_c}{2}\norm s^2.\label{eq:sweep-drop}
\end{align}

\emph{Gradient remaining after the sweep.}
When coordinate $i$ is updated, $\partial_i\Phi(z^{(i)})=0$. Each later update changes a coordinate $j>i$, so the term $-b_i/z_i$ in the $i$th derivative stays fixed. The remaining change comes entirely from the linear term $(Qz)_i$. Thus
\begin{align*}
 \partial_i\Phi(z^{(m)})
 &=\partial_i\Phi(z^{(i)})+\bigl[Q(z^{(m)}-z^{(i)})\bigr]_i\\
 &=\sum_{j>i}Q_{ij}s_j\\
 &=(Us)_i.
\end{align*}
Write $H_{\rm old}=H(z^{(0)})$ and $H_{\rm new}=H(z^{(m)})$. The gradient estimate and the decrease now give
\begin{align*}
 H_{\rm new}
 &\le \frac{\norm{Us}^2}{2\mu} && \text{by \eqref{eq:gap-gradient}}\\
 &\le \frac{M^2}{2\mu}\norm s^2 && \text{by }M=\norm U\\
 &\le \frac{M^2}{\mu\mu_c}D && \text{by \eqref{eq:sweep-drop}}\\
 &= \frac{M^2}{\mu\mu_c}(H_{\rm old}-H_{\rm new}).
\end{align*}
Multiplying by $\mu\mu_c$ and collecting the terms in $H_{\rm new}$ yields
\[
 (M^2+\mu\mu_c)H_{\rm new}\le M^2 H_{\rm old}.
\]
Division gives \eqref{eq:contraction}. Applying it to each successive sweep gives $H_k\le\rho^kH_0$. The distance estimate follows by taking $z=a^{(k)}$ in the lower bound for $H(z)$. When $M=0$, the first sweep ends with $\nabla\Phi(z^{(m)})=Us=0$, which characterizes the unique minimizer.
\end{proof}

\begin{corollary}[Finite residual stopping]\label{cor:stop}
Let $\ell=\min_i\ell_i>0$, where \eqref{eq:lower-coord} bounds each $a_i$ below on $\mathcal S$. With $b_{\max}=\max_i b_i$, set
\[
 L_\Phi=\norm Q+b_{\max}/\ell^2,\qquad
 K=\frac{R L_\Phi}{b_{\min}}\sqrt{\frac{2H_0}{\mu}}.
\]
If $0<\rho<1$, the residual after $k$ sweeps satisfies
\begin{equation}\label{eq:res-rate}
 \mathcal R(a^{(k)};Q,b)\le K\rho^{k/2}.
\end{equation}
For $K>\delta>0$, the test succeeds within
\begin{equation}\label{eq:sweep-count}
 k=\left\lceil\frac{2\log(K/\delta)}{\log(1/\rho)}\right\rceil
\end{equation}
sweeps. For $K\le\delta$ the initial iterate meets the tolerance.
\end{corollary}
\begin{proof}
The lower bound $a_i\ge\ell$ gives $\norm{\nabla^2\Phi(a)}\le\norm Q+b_{\max}/\ell^2=L_\Phi$ throughout $\mathcal S$. Set $h=a-a_\star$. Integrating from $a_\star$ and using $\nabla\Phi(a_\star)=0$ gives
\begin{align*}
 \norm{\nabla\Phi(a)}
 &=\left\|\int_0^1 \nabla^2\Phi(a_\star+th)h\,dt\right\|\\
 &\le \int_0^1 L_\Phi\norm h\,dt=L_\Phi\norm{a-a_\star}.
\end{align*}
The identity $e_i=(a_i/b_i)\partial_i\Phi(a)$ connects this bound to the observable residual. At $a=a^{(k)}$ it gives
\begin{align*}
 \mathcal R(a;Q,b)
 &\le (R/b_{\min})\norm{\nabla\Phi(a)}\\
 &\le (R L_\Phi/b_{\min})\norm{a-a_\star}\\
 &\le (R L_\Phi/b_{\min})\sqrt{2H_k/\mu}\\
 &\le K\rho^{k/2}.
\end{align*}
For $K>\delta>0$, taking logarithms gives
\[
 K\rho^{k/2}\le\delta
 \quad\Longleftrightarrow\quad
 (k/2)\log(1/\rho)\ge\log(K/\delta).
\]
The factor $\log(1/\rho)$ is positive, so rounding the lower bound on $k$ up gives \eqref{eq:sweep-count}. For $K\le\delta$, set $k=0$ in \eqref{eq:res-rate}.
\end{proof}

\subsection{Interpreting the sweep bound}
Algorithm 3 tests $\mathcal R$ directly. The constants $K$ and $\rho$ explain the cost of reaching a tolerance on the current subproblem. A diagonal $Q$ gives one-sweep convergence. For interacting tasks, the bound depends on the upper triangular part of $Q$ in the chosen order and on the curvature over the initial sublevel set. Those quantities can change at each training step, particularly when $\chi$ is small. Section~\ref{sec:research} counts that changing work across outer iterations.

\section{From an inner residual to outer convergence}\label{sec:outer}
Section~\ref{sec:solver} stops the inner calculation when every contribution is within a chosen tolerance of its budget. Here we show why those approximate weights suffice for training. The residual gives a lower bound on the improvement each task receives along the normalized direction. Smoothness turns that bound into a decrease for one step. Summing the decreases gives a rate of approach to Pareto stationarity.

\subsection{A certificate for simultaneous descent}
\begin{lemma}[Residual identities and alignment]\label{lem:alignment}
Let $Q=G^\top G$. Suppose $a>0$ satisfies $\mathcal R(a;Q,b)\le\delta<1$, with signed errors $e_i$ from \eqref{eq:residual}. Write
\[
 s=\one^\top a,\qquad w=a/s,\qquad v=Ga,\qquad d=Gw=v/s.
\]
The unnormalized direction satisfies the identities
\begin{equation}\label{eq:res-identities}
 g_i^\top v=\frac{b_i(1+e_i)}{a_i}>0,\qquad
 \norm v^2=\sum_i b_i(1+e_i).
\end{equation}
The normalized direction satisfies
\begin{equation}\label{eq:scaled-identities}
 g_i^\top d\ge\frac{(1-\delta)b_{\min}}{s^2},\qquad
 \norm d^2\le\frac{(1+\delta)B}{s^2}.
\end{equation}
In particular, with $c_\delta=(1-\delta)b_{\min}/[(1+\delta)B]$,
\begin{equation}\label{eq:align}
 g_i^\top d\ge c_\delta\norm d^2\quad(1\le i\le m),
 \qquad \chi(x)\le\norm d.
\end{equation}
\end{lemma}
\begin{proof}
The residual bound places $1+e_i$ between $1-\delta$ and $1+\delta$. Starting from the definition of $v$, we obtain
\begin{align*}
 g_i^\top v
 &=g_i^\top Ga=(Qa)_i\\
 &=\frac{b_i}{a_i}(1+e_i) && \text{by \eqref{eq:residual}}\\
 &\ge\frac{b_i}{a_i}(1-\delta)>0.
\end{align*}
Multiplying by $a_i$ and summing yields
\begin{align*}
 \norm v^2
 &=a^\top Qa\\
 &=\sum_i a_i(Qa)_i\\
 &=\sum_i b_i(1+e_i).
\end{align*}
Thus $(1-\delta)B\le\norm v^2\le(1+\delta)B$, since $\sum_i b_i=B$.

Now divide the direction by $s$. Since $s=\sum_j a_j\ge a_i$, the two quantities of interest satisfy
\begin{align*}
 g_i^\top d
 &=\frac{b_i(1+e_i)}{s a_i}
 \ge\frac{b_i(1-\delta)}{s a_i}
 \ge\frac{b_{\min}(1-\delta)}{s^2},\\
 \norm d^2&=\frac{\norm v^2}{s^2}
 \le\frac{(1+\delta)B}{s^2}.
\end{align*}
The same factor $1/s^2$ occurs in both bounds. Eliminating it yields
\begin{align*}
 g_i^\top d
 &\ge b_{\min}(1-\delta)\,\frac{1}{s^2}\\
 &\ge \frac{b_{\min}(1-\delta)}{(1+\delta)B}\norm d^2
 = c_\delta\norm d^2.
\end{align*}
Finally, $w=a/s$ belongs to the simplex. Evaluating the minimum at this feasible point gives $\chi(x)=\min_{u\in\simplex}\norm{Gu}\le\norm{Gw}=\norm d$.
\end{proof}

The constant $c_\delta$ depends only on the budgets and the chosen tolerance. The same positive tolerance therefore supplies the same alignment bound at every training step.

\begin{assumption}\label{ass:smooth}
For each task $i$, the objective $f_i$ is continuously differentiable on $\R^p$ and admits a finite lower bound $f_i^{\rm low}$. There is a common $L>0$ such that, for every $x,y\in\R^p$,
\[
 \norm{\nabla f_i(x)-\nabla f_i(y)}\le L\norm{x-y}.
\]
When the tasks have individual constants $L_i$, one can use $L=\max_i L_i$ if this maximum is positive.
\end{assumption}

\begin{procedure}{Algorithm 4. Normalized Nash descent with a relative residual tolerance}
\begin{enumerate}[leftmargin=*]
\item Choose the starting point $x_0$, budgets $b>0$, tolerance $0<\delta<1$, and step size $\eta=c_\delta/L$.
\item At $x_t$, compute the task-gradient matrix $G_t$ and $Q_t=G_t^\top G_t$.
\item If $\chi(x_t)=0$, return $x_t$. Otherwise run Algorithm 3 until it returns $a_t>0$ with $\mathcal R(a_t;Q_t,b)\le\delta$.
\item Use the normalized weights to take the step
\[
 w_t=\frac{a_t}{\one^\top a_t},\qquad d_t=G_tw_t,\qquad x_{t+1}=x_t-\eta d_t.
\]
Repeat from Step 2.
\end{enumerate}
\end{procedure}

For a numerical stationarity test, any feasible simplex vector $w$ with $\norm{G_tw}\le\varepsilon$ certifies $\chi(x_t)\le\varepsilon$. One can obtain such a vector by minimizing the simplex quadratic $\norm{G_tw}^2$. After a successful inner solve, $\norm{d_t}\le\varepsilon$ provides the same certificate. If a sweep cap is reached before meeting either stopping criterion, report the inner solve as uncertified.

\subsection{From decrease of each objective to Pareto stationarity}
\begin{theorem}[Deterministic nonconvex convergence]\label{thm:outer}
Under Assumption~\ref{ass:smooth}, Algorithm 4 is well-defined at every nonstationary iteration. Every performed update decreases each task's objective according to
\begin{align}
 f_i(x_{t+1})
 &\le f_i(x_t)-\frac{c_\delta^2}{2L}\norm{d_t}^2\nonumber\\
 &\le f_i(x_t)-\frac{c_\delta^2}{2L}\chi(x_t)^2.\label{eq:descent}
\end{align}
Write $\Delta=\min_i[f_i(x_0)-f_i^{\rm low}]$. If $T\ge1$ updates are performed, then
\begin{equation}\label{eq:outer-rate}
 \min_{0\le t<T}\chi(x_t)^2\le\frac{2L\Delta}{c_\delta^2T}.
\end{equation}
With equal budgets, this becomes
\begin{equation}\label{eq:equal-rate}
 \min_{0\le t<T}\chi(x_t)^2
 \le\frac{2Lm^2(1+\delta)^2\Delta}{(1-\delta)^2 T}.
\end{equation}
Along an infinite run, both $\norm{d_t}$ and $\chi(x_t)$ tend to zero. Every accumulation point is Pareto stationary.
\end{theorem}
\begin{proof}
At any point with $\chi(x_t)>0$, Theorem~\ref{thm:exist} provides a positive exact solution. Corollary~\ref{cor:stop} and the diagonal case of Theorem~\ref{thm:ccd} ensure that Algorithm 3 reaches the prescribed tolerance in finitely many sweeps.

\emph{A bound for one step.}
Consider any displacement $h$ from $x$. Subtract the first-order term and apply the fundamental theorem of calculus. Cauchy--Schwarz and the gradient Lipschitz bound give
\begin{align*}
 &f_i(x+h)-f_i(x)-\nabla f_i(x)^\top h\\
 &=\int_0^1[\nabla f_i(x+\tau h)-\nabla f_i(x)]^\top h\,d\tau\\
 &\le\int_0^1\norm{\nabla f_i(x+\tau h)-\nabla f_i(x)}\,\norm h\,d\tau\\
 &\le\int_0^1 L\tau\norm h^2\,d\tau
 =\tfrac L2\norm h^2.
\end{align*}
Set $g_{i,t}=\nabla f_i(x_t)$ and $h=-\eta d_t$. The alignment inequality \eqref{eq:align} bounds the linear part of the step, giving
\begin{align*}
 f_i(x_{t+1})-f_i(x_t)
 &\le-\eta g_{i,t}^\top d_t+\tfrac{L\eta^2}{2}\norm{d_t}^2\\
 &\le-\left(\eta c_\delta-\tfrac{L\eta^2}{2}\right)\norm{d_t}^2\\
 &=-\frac{c_\delta^2}{2L}\norm{d_t}^2
 \qquad (\eta=c_\delta/L).
\end{align*}
Since $\chi(x_t)\le\norm{d_t}$, this proves both inequalities in \eqref{eq:descent}.

\emph{Summing over iterations.}
Fix any task $i$. Summing the one-step inequality and applying the lower bound yields
\begin{align*}
 \frac{c_\delta^2}{2L}\sum_{t=0}^{T-1}\norm{d_t}^2
 &\le\sum_{t=0}^{T-1}\bigl[f_i(x_t)-f_i(x_{t+1})\bigr]\\
 &= f_i(x_0)-f_i(x_T)\\
 &\le f_i(x_0)-f_i^{\rm low}.
\end{align*}
The left side is common to all tasks. Taking the smallest right side and then using the minimum-average inequality gives
\begin{align*}
 \min_{0\le t<T}\chi(x_t)^2
 &\le \frac1T\sum_{t=0}^{T-1}\chi(x_t)^2\\
 &\le \frac1T\sum_{t=0}^{T-1}\norm{d_t}^2\\
 &\le\frac{2L\Delta}{c_\delta^2 T}.
\end{align*}
For equal budgets $b_i=b_0>0$, use $b_{\min}=b_0$ and $B=mb_0$ to obtain
\[
 c_\delta=\frac{1-\delta}{m(1+\delta)},\qquad
 \frac1{c_\delta^2}=\frac{m^2(1+\delta)^2}{(1-\delta)^2}.
\]
Substitution gives \eqref{eq:equal-rate}.

\emph{Infinite runs and accumulation points.}
The partial sums of the nonnegative series $\sum_t\norm{d_t}^2$ are bounded by $2L\Delta/c_\delta^2$. They therefore converge, and their consecutive differences tend to zero. It follows that
\[
 \norm{d_t}^2\longrightarrow0,\qquad
 0\le\chi(x_t)\le\norm{d_t}\longrightarrow0.
\]
To pass to a limit point, choose a minimizer $w_y$ in the definition of $\chi(y)$. Such a minimizer exists by compactness of the simplex. Then
\begin{align*}
 \chi(x)
 &\le\norm{G(x)w_y}\\
 &\le\norm{G(y)w_y}
       +\sum_i(w_y)_i\norm{\nabla f_i(x)-\nabla f_i(y)}\\
 &\le\chi(y)+L\norm{x-y}\sum_i(w_y)_i\\
 &=\chi(y)+L\norm{x-y}.
\end{align*}
Interchanging $x$ and $y$ gives $|\chi(x)-\chi(y)|\le L\norm{x-y}$. Thus if $x_{t_k}\to\bar x$, continuity gives $\chi(\bar x)=\lim_k\chi(x_{t_k})=0$, which is Pareto stationarity.
\end{proof}

For the simplex-normalized step in Algorithm 4, the rate gives $O(\varepsilon^{-2})$ iterations to reach $\chi\le\varepsilon$. When accuracy is expressed as $\chi^2\le\varepsilon$, the same bound gives $O(\varepsilon^{-1})$ iterations.

\subsection{Choosing steps without a known smoothness constant}
Backtracking replaces the fixed step $\eta=c_\delta/L$ by a test of the actual objective decreases. Choose $0<\sigma<1$, an initial trial $\eta_{\max}>0$, and a reduction factor $0<\theta<1$. At each iteration, take the first $\eta$ in the sequence $\eta_{\max},\theta\eta_{\max},\ldots$ satisfying
\begin{equation}\label{eq:armijo}
 f_i(x-\eta d)\le f_i(x)-\sigma\eta g_i^\top d
 \qquad\text{for every }i.
\end{equation}
The smoothness estimate above bounds the left side. After reserving the fraction $\sigma$ of the linear decrease required by the test, we obtain
\begin{align*}
 f_i(x-\eta d)-f_i(x)+\sigma\eta g_i^\top d
 &\le -(1-\sigma)\eta g_i^\top d+\tfrac{L\eta^2}{2}\norm d^2\\
 &\le-\eta\left((1-\sigma)c_\delta-\tfrac{L\eta}{2}\right)\norm d^2.
\end{align*}
This is nonpositive whenever $\eta\le\eta_{\rm suff}:=2(1-\sigma)c_\delta/L$. Since the trial steps decrease geometrically to zero, one is accepted after finitely many reductions. If the first trial succeeds, the accepted step equals $\eta_{\max}$. Otherwise the previous trial failed, so $\eta/\theta>\eta_{\rm suff}$. In both cases,
\[
 \eta\ge\eta_{\min}:=\min\{\eta_{\max},\theta\eta_{\rm suff}\}>0.
\]
For any task, the accepted step therefore gives
\begin{align*}
 f_i(x_t)-f_i(x_{t+1})
 &\ge\sigma\eta\,g_{i,t}^\top d_t && \text{by \eqref{eq:armijo}}\\
 &\ge\sigma\eta_{\min}c_\delta\norm{d_t}^2 && \text{by \eqref{eq:align}}\\
 &\ge\sigma\eta_{\min}c_\delta\chi(x_t)^2.
\end{align*}
Summing over the $T$ updates and using the definition of $\Delta$ yields
\[
 \min_{0\le t<T}\chi(x_t)^2
 \le\frac1T\sum_{t=0}^{T-1}\chi(x_t)^2
 \le\frac{\Delta}{\sigma\eta_{\min}c_\delta T}.
\]

\section{Changing matrices, regularization, and noisy gradients}\label{sec:perturb}
Section~\ref{sec:outer} certifies descent using the gradients at the current iterate. In practice, weights may come from an earlier iterate, a regularized matrix, or a batch of approximate gradients. We first compute the change in the residual when only the matrix used to find the weights changes. We then consider errors in the direction actually taken and trace their effect through the outer convergence proof.

\subsection{A residual certificate for matrix error and stale weights}
Fix the weights $a$. Write $\widehat Q$ for the matrix used to compute them and $Q$ for the matrix against which we want to certify descent.

\begin{proposition}[Residual transfer]\label{prop:matrix-error}
Let $\widehat Q=Q+E$ with $Q$ and $\widehat Q$ symmetric. If $a>0$, $b>0$, and $\mathcal R(a;\widehat Q,b)\le\widehat\delta$, then
\begin{equation}\label{eq:res-transfer}
 \mathcal R(a;Q,b)
 \le\widehat\delta+
 \frac{\norm a_\infty\norm a}{b_{\min}}\norm E.
\end{equation}
When $Q=G^\top G$ and this upper bound is less than one, $g_i^\top Ga>0$ for every task $i$.
\end{proposition}
\begin{proof}
Substituting $Q=\widehat Q-E$ into each coordinate gives
\begin{align*}
 e_i(a;Q,b)
 &= \frac{a_i(Qa)_i}{b_i}-1\\
 &= \frac{a_i(\widehat Qa)_i}{b_i}-1-\frac{a_i(Ea)_i}{b_i}\\
 &= e_i(a;\widehat Q,b)-\frac{a_i(Ea)_i}{b_i}.
\end{align*}
Since $|(Ea)_i|\le\norm{Ea}$, the triangle inequality gives
\begin{align*}
 |e_i(a;Q,b)|
 &\le |e_i(a;\widehat Q,b)|+\frac{a_i}{b_i}|(Ea)_i|\\
 &\le \widehat\delta+\frac{\norm a_\infty}{b_{\min}}\norm{Ea}\\
 &\le \widehat\delta+\frac{\norm a_\infty\norm a}{b_{\min}}\norm E.
\end{align*}
Taking the maximum over $i$ proves \eqref{eq:res-transfer}. If the resulting bound is less than one, then $1+e_i(a;Q,b)>0$. Rearranging the defining equation for $e_i$ gives
\[
 g_i^\top Ga=(Qa)_i=\frac{b_i}{a_i}\bigl(1+e_i(a;Q,b)\bigr)>0.
 \qedhere
\]
\end{proof}

\paragraph{Reusing weights at a new iterate.}
Suppose $a$ is certified at $x$ and we want to use it at $y$. In \eqref{eq:res-transfer} take $\widehat Q=Q(x)$ and $Q=Q(y)$. If $Q(y)$ is already available, a single product $Q(y)a$ gives the current residual exactly. A bound on the change in $Q$ gives another way to certify reuse.

Assume that every $g_i$ is $L$-Lipschitz and that $\norm{g_i(z)}\le M_g$ for $z=x,y$. Put $D=G(y)-G(x)$. The columns of $D$ are the changes in the task gradients. Bounding the operator norm by the Frobenius norm gives
\begin{align*}
 \norm{G(z)}^2
 &\le \sum_i\norm{g_i(z)}^2
 \le mM_g^2,\\
 \norm D^2
 &\le \sum_i\norm{g_i(y)-g_i(x)}^2
 \le mL^2\norm{y-x}^2.
\end{align*}
Adding and subtracting $G(y)^\top G(x)$ gives
\begin{align*}
 Q(y)-Q(x)
 &=G(y)^\top[G(y)-G(x)]+[G(y)-G(x)]^\top G(x)\\
 &=G(y)^\top D+D^\top G(x),\\[2pt]
 \norm{Q(y)-Q(x)}
 &\le \bigl(\norm{G(y)}+\norm{G(x)}\bigr)\norm D\\
 &\le 2mM_gL\norm{y-x}.
\end{align*}
Thus weights with residual at most $\widehat\delta$ at $x$ remain within a prescribed tolerance $\delta$ at $y$ whenever
\[
 \widehat\delta+
 \frac{2mM_gL\,\norm a_\infty\norm a}{b_{\min}}\,\norm{y-x}
 \le\delta.
\]
Using this test before forming $Q(y)$ requires bounds on $M_g$ and $L$.

\subsection{Sensitivity of the exact weights}
The previous bound holds the weights fixed and measures the change in their residual. We now allow the weights to adjust and bound the distance between the two exact minimizers.

\begin{proposition}[Local perturbation bound]\label{prop:sensitivity}
Let $Q,Q'\succeq0$ both satisfy the nonnegative-kernel condition of Theorem~\ref{thm:map}. Let $a,a'$ be the respective minimizers for a common $b>0$. Choose $R$ such that every coordinate of both minimizers is at most $R$ and put $\mu=b_{\min}/R^2$. Then
\begin{equation}\label{eq:sensitivity}
 \norm{a'-a}\le\frac{\norm{Q'-Q}\norm a}{\mu}.
\end{equation}
\end{proposition}
\begin{proof}
Put $h=a'-a$ and $z_\tau=a+\tau h$ for $0\le\tau\le1$. Each coordinate of $z_\tau$ lies in $(0,R]$. The curvature along the segment therefore satisfies
\begin{align*}
 h^\top\nabla^2\Phi_{Q',b}(z_\tau)h
 &= h^\top Q'h+\sum_i\frac{b_i}{z_{\tau,i}^2}h_i^2\\
 &\ge \frac{b_{\min}}{R^2}\sum_i h_i^2
 =\mu\norm h^2 .
\end{align*}
Each endpoint minimizes its own objective. Evaluating both gradients under the new matrix $Q'$ gives
\begin{align*}
 \nabla\Phi_{Q',b}(a') &= Q'a'-b/a'=0,\\
 \nabla\Phi_{Q',b}(a)  &= (Q'-Q)a+\underbrace{Qa-b/a}_{=\,0}
                        =(Q'-Q)a .
\end{align*}
Integrating along the segment and substituting the endpoint gradients gives
\begin{align*}
 \mu\norm h^2
 &\le\int_0^1 h^\top\nabla^2\Phi_{Q',b}(z_\tau)h\,d\tau\\
 &=h^\top\bigl[\nabla\Phi_{Q',b}(a')-\nabla\Phi_{Q',b}(a)\bigr]\\
 &=-h^\top(Q'-Q)a\\
 &\le\norm h\,\norm{Q'-Q}\,\norm a .
\end{align*}
For $h\ne0$, division by $\mu\norm h$ yields \eqref{eq:sensitivity}. For $h=0$, its left side is zero.
\end{proof}
The factor $\norm a/\mu$ describes the accuracy of a warm start. For the same change in $Q$, the bound allows a larger change in the optimal weights when the weights are large or the curvature is small.

\subsection{Ridge regularization changes the certificate}
Adding $\lambda I$, with $\lambda>0$, makes $Q_\lambda=Q+\lambda I$ positive definite. By Theorem~\ref{thm:exist} the regularized problem has a unique positive solution, including at points where $\chi=0$. For each task, the regularized equation expands to
\begin{align*}
 (Q_\lambda a)_i &= b_i/a_i,\\
 (Qa)_i+\lambda a_i &= b_i/a_i,\\
 a_i(Qa)_i &= b_i-\lambda a_i^2 .
\end{align*}
We can now read off the residual and directional utility for the original matrix,
\begin{equation}\label{eq:ridge}
 \begin{aligned}
 e_i(a;Q,b) &= \frac{b_i-\lambda a_i^2}{b_i}-1
              =-\frac{\lambda a_i^2}{b_i},\\
 g_i^\top Ga &= (Qa)_i=\frac{b_i}{a_i}-\lambda a_i .
 \end{aligned}
\end{equation}
In particular $g_i^\top Ga>0$ holds exactly when $\lambda a_i^2<b_i$.

For an approximate regularized solution, let $\mathcal R(a;Q_\lambda,b)\le\widehat\delta$. Subtracting the ridge term yields a bound on the original residual,
\begin{align*}
 e_i(a;Q,b)
 &=e_i(a;Q_\lambda,b)-\frac{\lambda a_i^2}{b_i},\\
 |e_i(a;Q,b)|
 &\le\widehat\delta+\frac{\lambda a_i^2}{b_i},\\
 \mathcal R(a;Q,b)
 &\le\widehat\delta+\lambda\max_i\frac{a_i^2}{b_i}.
\end{align*}
Testing the original residual tells us whether to accept the step or reduce $\lambda$ and solve again.

\subsection{A deterministic gradient-error theorem}
Approximate task gradients affect both the weights and the update direction. At iteration $t$, write $g_{i,t}=\nabla f_i(x_t)$ and let $\widehat g_{i,t}$ denote its approximation. Assume that the realized errors satisfy
\begin{equation}\label{eq:grad-error}
 \max_i\norm{\widehat g_{i,t}-g_{i,t}}\le\zeta_t.
\end{equation}
Form $\widehat G_t$ from the approximate gradients and set $\widehat Q_t=\widehat G_t^\top\widehat G_t$. Given a positive $a_t$ with $\mathcal R(a_t;\widehat Q_t,b)\le\delta<1$, the update is
\begin{align*}
 w_t &= a_t/(\one^\top a_t),\\
 \widehat d_t &= \widehat G_t w_t,\\
 x_{t+1} &= x_t-\eta\widehat d_t .
\end{align*}

\begin{theorem}[Error-dependent stationarity bound]\label{thm:noisy}
Suppose Assumption~\ref{ass:smooth} and \eqref{eq:grad-error} hold, and a certified positive $a_t$ is obtained at every iteration. With $c=c_\delta$ and $\eta=c/L$, each objective satisfies
\begin{equation}\label{eq:noisy-descent}
 f_i(x_{t+1})\le f_i(x_t)-\frac{\eta c}{4}\norm{\widehat d_t}^2
 +\frac\eta c\zeta_t^2.
\end{equation}
For $T\ge1$ and $\Delta=\min_i[f_i(x_0)-f_i^{\rm low}]$, the true gradients satisfy
\begin{equation}\label{eq:noisy-rate}
 \min_{0\le t<T}\chi(x_t)^2
 \le\frac{8L\Delta}{c^2T}
 +\left(\frac8{c^2}+2\right)\frac1T\sum_{t<T}\zeta_t^2 .
\end{equation}
\end{theorem}
\begin{proof}
\emph{1. Bound the true directional utility.}
Write $z_t=\norm{\widehat d_t}$. Lemma~\ref{lem:alignment} gives $\widehat g_{i,t}^\top\widehat d_t\ge c z_t^2$ for the gradients used in the inner solve. Adding back the gradient error gives
\begin{align*}
 g_{i,t}^\top\widehat d_t
 &=\widehat g_{i,t}^\top\widehat d_t+
       (g_{i,t}-\widehat g_{i,t})^\top\widehat d_t\\
 &\ge c z_t^2-\norm{g_{i,t}-\widehat g_{i,t}}\,z_t\\
 &\ge c z_t^2-\zeta_t z_t .
\end{align*}

\emph{2. Bound the change in one objective.}
With the preceding estimate, the smoothness inequality from Theorem~\ref{thm:outer} gives
\begin{align*}
 f_i(x_{t+1})-f_i(x_t)
 &\le-\eta g_{i,t}^\top\widehat d_t+\tfrac12 L\eta^2 z_t^2\\
 &\le-\eta(c-\tfrac12 L\eta)z_t^2+\eta\zeta_t z_t\\
 &=-\tfrac12\eta c z_t^2+\eta\zeta_t z_t
       \qquad (L\eta=c).
\end{align*}
To bound the product of the gradient error and the step length, use
\[
 0\le\left(\frac{\sqrt c}{2}z_t-\frac{\zeta_t}{\sqrt c}\right)^2
   =\frac c4 z_t^2-\zeta_t z_t+\frac{\zeta_t^2}{c}.
\]
Replacing $\zeta_t z_t$ by its upper bound $c z_t^2/4+\zeta_t^2/c$ proves \eqref{eq:noisy-descent}.

\emph{3. Add the decreases over time.}
For each task, summing \eqref{eq:noisy-descent} and rearranging yields
\begin{align*}
 \frac{\eta c}{4}\sum_{t<T}z_t^2
 &\le \sum_{t<T}\bigl[f_i(x_t)-f_i(x_{t+1})\bigr]
             +\frac{\eta}{c}\sum_{t<T}\zeta_t^2\\
 &= f_i(x_0)-f_i(x_T)+\frac{\eta}{c}\sum_{t<T}\zeta_t^2\\
 &\le f_i(x_0)-f_i^{\rm low}+\frac{\eta}{c}\sum_{t<T}\zeta_t^2 .
\end{align*}
Minimizing over tasks on the right and dividing by $\eta c/4$ gives
\begin{equation}\label{eq:noisy-sum}
 \sum_{t<T}z_t^2
 \le\frac{4\Delta}{\eta c}+\frac4{c^2}\sum_{t<T}\zeta_t^2.
\end{equation}

\emph{4. Return to the true stationarity measure.}
The simplex weights are nonnegative and sum to one. Consequently
\begin{align*}
 \chi(x_t)
 &\le\norm{G_tw_t}\\
 &\le\norm{\widehat G_t w_t}
       +\left\|\sum_i w_{i,t}(g_{i,t}-\widehat g_{i,t})\right\|\\
 &\le z_t+\sum_i w_{i,t}\zeta_t\\
 &=z_t+\zeta_t.
\end{align*}
Applying $(z_t+\zeta_t)^2\le2z_t^2+2\zeta_t^2$ and then \eqref{eq:noisy-sum} gives
\begin{align*}
 \min_{t<T}\chi(x_t)^2
 &\le \frac1T\sum_{t<T}\chi(x_t)^2\\
 &\le \frac2T\sum_{t<T}z_t^2+\frac2T\sum_{t<T}\zeta_t^2\\
 &\le \frac{8\Delta}{\eta cT}
        +\left(\frac8{c^2}+2\right)\frac1T\sum_{t<T}\zeta_t^2 .
\end{align*}
Setting $\eta=c/L$ gives \eqref{eq:noisy-rate}.
\end{proof}
The proof uses the realized gradient errors and allows the weights and update to use the same batch. For a random sequence of batches, \eqref{eq:noisy-rate} holds on every event where both the error bounds and the inner certificates hold. If these conditions hold almost surely and the squared errors are integrable, taking expectations on both sides yields an expected bound. Applying the theorem to a particular sampling scheme therefore requires estimates of the squared errors and a rule for sampled gradients that admit no strictly positive solution.

\section{From the correspondence to a research study}\label{sec:research}

A solver can take fewer sweeps per update and still have little effect on training time if computing task gradients accounts for most of the work. We combine the two iteration bounds to see which costs the theory controls. This calculation also identifies what to measure in a comparison of solvers.

\subsection{The proposed contribution and its relation to prior results}

The quadratic logarithmic objective appears in both risk budgeting and multi-task learning \citep{spinu2013,li2026toppo}. Normalized Nash aggregation also admits a convergence analysis within the framework of \citet{hu2026}. The question pursued here concerns the accuracy of the weights \emph{at the moment the solver stops}. Corollary~\ref{cor:stop} bounds the sweeps needed to reach an observable relative budget error. Theorem~\ref{thm:outer} takes that same error as input and bounds the number of training steps. Together they allow a comparison at a common level of accuracy, even when two solvers use different internal iterations.

\subsection{Outer iterations and total arithmetic work}

Suppose $m$ tasks share $p$ parameters. At iteration $t$, let $k_t$ denote the number of coordinate sweeps. Once the $m$ task gradients are available, the main operations have the following costs.
\begin{center}
\begin{tabular}{@{}ll@{}}
\toprule
Operation at iteration $t$ & Dense arithmetic cost\\
\midrule
Form the $m\times m$ Gram matrix $G_t^\top G_t$ & $O(pm^2)$\\
Run $k_t$ complete sweeps and check the residual & $O((1+k_t)m^2)$\\
Form the update direction $G_tw_t$ & $O(pm)$\\
Run a separate simplex stationarity test, if used & $C_{{\rm stat},t}$\\
\bottomrule
\end{tabular}
\end{center}
The Gram matrix uses $O(m^2)$ storage. Keeping all gradients uses another $O(pm)$. Checking the norm of an already computed direction costs $O(p)$, which is included in the cost of forming it.

Summing the table over $T$ updates gives
\begin{equation}\label{eq:work}
 O\!\left(Tpm^2 + m^2\sum_{t=0}^{T-1}k_t + Tpm\right)
 +\sum_{t=0}^{T-1} C_{{\rm stat},t} .
\end{equation}
Here the initial residual checks cost $O(Tm^2)$, already covered by $O(Tpm^2)$ for $p\ge1$. In addition, the algorithm evaluates $mT$ task gradients. A warm start can reduce the sum of the sweep counts. Forming the Gram matrix and computing the gradients still contribute at each update.

\begin{corollary}[A conditional total-work bound]\label{cor:total}
Fix a residual tolerance $0<\delta<1$ and a stationarity tolerance $\varepsilon>0$. For every visited nondiagonal, nonstationary subproblem, suppose the constants of Corollary~\ref{cor:stop} satisfy
\[
 K_t\le\overline K<\infty,\qquad
 \rho_t\le\overline\rho<1,\qquad
 \overline K>0,\quad\overline\rho>0 .
\]
Diagonal problems take at most one sweep. Define
\begin{align*}
 \overline k
 &=\max\!\left\{1,\left\lceil
       \frac{2\max\{0,\log(\overline K/\delta)\}}
            {\log(1/\overline\rho)}
     \right\rceil\right\},\\[2pt]
 T_\varepsilon
 &=\max\!\left\{1,\left\lceil
       \frac{2L\Delta}{c_\delta^2\varepsilon^2}
     \right\rceil\right\}.
\end{align*}
Each inner solve uses at most $\overline k$ sweeps. Under Assumption~\ref{ass:smooth}, Algorithm 4 visits an iterate with $\chi(x_t)\le\varepsilon$ among the first $T_\varepsilon$ evaluated iterates, including any earlier termination at $\chi=0$. For $T$ performed updates, the sweep count in \eqref{eq:work} is at most $T\overline k$.
\end{corollary}

\begin{proof}
\emph{Sweeps per update.}
Consider a nondiagonal subproblem with $K_t>\delta$. The uniform bounds give
\begin{align*}
 \log(K_t/\delta)&\le\log(\overline K/\delta),\\
 \log(1/\rho_t)&\ge\log(1/\overline\rho)>0 .
\end{align*}
The first line bounds the numerator of the sweep estimate. The second bounds the denominator away from zero. Substitution in Corollary~\ref{cor:stop} yields
\begin{align*}
 k_t
 &\le
 \left\lceil\frac{2\log(K_t/\delta)}{\log(1/\rho_t)}\right\rceil\\
 &\le
 \left\lceil\frac{2\log(\overline K/\delta)}
                   {\log(1/\overline\rho)}\right\rceil
 \le\overline k .
\end{align*}
For $K_t\le\delta$, the initial vector already passes the test. The term $1$ in the definition of $\overline k$ covers a diagonal subproblem.

\emph{Training steps.}
If the run terminates earlier at $\chi=0$, the stationarity target holds there. Otherwise Theorem~\ref{thm:outer}, applied with $T=T_\varepsilon$, gives
\begin{align*}
 \min_{0\le t<T_\varepsilon}\chi(x_t)^2
 &\le \frac{2L\Delta}{c_\delta^2 T_\varepsilon}\\
 &\le \varepsilon^2
 \qquad\text{since }c_\delta^2 T_\varepsilon\varepsilon^2
       \ge 2L\Delta .
\end{align*}
Taking the nonnegative square root gives the target. Finally,
\[
 \sum_{t=0}^{T-1} k_t
 \le \sum_{t=0}^{T-1}\overline k
 = T\overline k .
 \qedhere
\]
\end{proof}

The uniform bounds on $K_t$ and $\rho_t$ are additional assumptions about the inner problems. Smooth task objectives can have nearly opposing gradients. Appendix~\ref{app:examples} shows explicitly how the canonical weights grow in this regime.

\section{Conclusion}\label{sec:conclusion}

In Nash-MTL, the product of a task's weight and its directional utility equals its budget. In a risk-budgeting portfolio, the same product gives an asset's contribution to variance, up to a common scale. Once we match the matrices, normalizing the Nash weights gives the portfolio. Fixing the total quadratic risk recovers the Nash weights. A positive solution exists exactly when the gradients admit a direction of strict improvement for every task.

The link matters during training because a coordinate solver can measure each budget error as it runs. Once all relative errors fall below a fixed tolerance smaller than one, the normalized gradient mixture has a quantitative descent guarantee for every task. The remaining computational question is how many sweeps this requires across an entire run. Nearly opposing gradients provide a concrete test case. Their normalized mixture tends to zero while the canonical weights grow without bound, so the practical value of a warm start depends on how much the gradient geometry changes between updates.

\section{Disclosure}
\label{sec:disclosure}
The proof strategy is produced by OpenAI’s GPT-6.2 Astra Ultra through Codex in
response to prompts from the author. Codex was also used to revise
the exposition and prepare the LaTeX manuscript. The author selected the
problem, directed the interactions and revisions, and is the sole named author. The AI system is acknowledged as a reasoning and writing tool, not
as an author. This disclosure is not a substitute for independent expert
mathematical review.
\bibliographystyle{plainnat}
\begingroup
\small
\raggedright
\setlength{\bibsep}{4pt plus 1pt minus 1pt}
\bibliography{references}
\endgroup

\clearpage
\appendix
\section{Worked examples and boundary cases}\label{app:examples}
The first example follows a pair of weights from bargaining to portfolio allocation and back. The remaining examples vary the angle between two gradients to show when the weights exist and what happens as a common descent direction disappears.

\subsection{Orthogonal gradients and the two weight scales}
Take $g_1=(1,0)^\top$, $g_2=(0,2)^\top$, and $b=(1,1)^\top$. The two gradient columns give
\[
 G=\begin{pmatrix}1&0\\0&2\end{pmatrix},\qquad
 Q=G^\top G=\begin{pmatrix}1&0\\0&4\end{pmatrix}.
\]
Thus each coordinate of $Qa=b/a$ can be solved independently. Taking the positive square roots gives
\begin{align*}
 a_1=1/a_1 &\quad\Longrightarrow\quad a_1^2=1
             \quad\Longrightarrow\quad a_1=1,\\
 4a_2=1/a_2 &\quad\Longrightarrow\quad a_2^2=1/4
             \quad\Longrightarrow\quad a_2=1/2 .
\end{align*}
For these weights, the two directions and the scale relating them are
\begin{align*}
 v &= Ga=(1,1)^\top, & \norm v^2&=1^2+1^2=2=B,\\
 s &= a_1+a_2=3/2, & w&=a/s=(2/3,1/3)^\top,\\
 d &= Gw=v/s=(2/3,2/3)^\top .
\end{align*}
For task 2, twice the directional utility is balanced by half the weight. This gives equal weighted gains for the two tasks, since
\begin{align*}
 g_1^\top v&=1, & a_1 g_1^\top v&=1\cdot1=1,\\
 g_2^\top v&=2, & a_2 g_2^\top v&=\tfrac12\cdot2=1 .
\end{align*}
To see the same balance in the portfolio, first compute $Qw=(2/3,4/3)^\top$. Multiplying by the corresponding holdings gives
\begin{align*}
 w_1(Qw)_1 &= (2/3)(2/3)=4/9,\\
 w_2(Qw)_2 &= (1/3)(4/3)=4/9,\\
 w^\top Qw &= 4/9+4/9=8/9 .
\end{align*}
Each asset contributes the fraction $(4/9)/(8/9)=1/2$. Applying the inverse map to this portfolio gives
\begin{align*}
 \sqrt{\frac{B}{w^\top Qw}}\,w
 &=\sqrt{\frac{2}{8/9}}\,(2/3,1/3)^\top\\
 &=\tfrac32(2/3,1/3)^\top=(1,1/2)^\top.
\end{align*}

\subsection{Identical gradients and a singular Gram matrix}
Let $g_1=g_2=u$ with $\norm u=1$ and $b=\one$. Then $Q$ is the matrix of ones and both coordinates of $Qa$ equal $a_1+a_2$. The two budget equations give
\begin{align*}
 1/a_1 &= a_1+a_2=1/a_2,\\
 a_1 &= a_2=:\alpha,\\
 2\alpha &=1/\alpha,\qquad \alpha=1/\sqrt2 .
\end{align*}
The canonical direction is $v=2\alpha u=\sqrt2\,u$. For \emph{every} simplex vector $w$, including the solution $w=(1/2,1/2)^\top$,
\[
 Gw=(w_1+w_2)u=u,\qquad \chi=\norm u=1.
\]
The Gram matrix has rank one. Its null vectors are the multiples of $(1,-1)^\top$, and the only nonnegative vector among them is zero. This is precisely the existence condition in Theorem~\ref{thm:exist}.

\subsection{Exactly opposing gradients}
Let $g_1=u$, $g_2=-u$, $\norm u=1$, and $b=\one$. Now
\[
 Q=\begin{pmatrix}1&-1\\-1&1\end{pmatrix},
 \qquad G(1/2,1/2)^\top=\tfrac12u-\tfrac12u=0 .
\]
It follows that $\chi=0$. To see directly what fails in the weight problem, take $a=(t,t)^\top$, $t>0$. Along this ray,
\begin{align*}
 Qa &= (t-t,-t+t)^\top=0,\\
 \Phi_{Q,\one}(a)
 &=\tfrac12a^\top Qa-\log t-\log t\\
 &=-2\log t\longrightarrow-\infty\qquad(t\to\infty).
\end{align*}
Theorem~\ref{thm:exist} therefore rules out a positive solution to $Qa=\one/a$.

For comparison, a ridge $\lambda>0$ makes the regularized objective strictly convex with a unique minimizer. Swapping the coordinates leaves that objective unchanged, so uniqueness gives $a=(\alpha,\alpha)^\top$. Since $Qa=0$, either coordinate equation reduces to
\begin{align*}
 \lambda\alpha &=1/\alpha,\\
 \alpha &=1/\sqrt\lambda,\\
 Ga &=\alpha(u-u)=0,\\
 e_i(a;Q,\one) &=a_i(Qa)_i-1=-1 .
\end{align*}
The ridge solver thus returns finite positive weights at a point where every task's directional utility is zero.

\subsection{Nearly opposing gradients}
Let $0<\theta<\pi/2$. Define $g_1=(1,0)^\top$, $g_2=(-\cos\theta,\sin\theta)^\top$, and $c=\cos\theta$. With $b=\one$,
\[
 Q=\begin{pmatrix}1&-c\\-c&1\end{pmatrix}.
\]
The eigenvalues $1-c$ and $1+c$ are positive. Symmetry and uniqueness again give $a_1=a_2=\alpha$. This time the budget equation gives
\begin{align*}
 (1-c)\alpha &=1/\alpha,\\
 \alpha &=1/\sqrt{1-c},\\
 w &= (\alpha,\alpha)^\top/(2\alpha)=(1/2,1/2)^\top,\\
 d &=Gw=\tfrac12(1-c,\sin\theta)^\top.
\end{align*}
Its squared length is
\begin{align*}
 \norm d^2
 &=\tfrac14\bigl((1-c)^2+\sin^2\theta\bigr)\\
 &=\tfrac14\bigl(1-2c+c^2+1-c^2\bigr)\\
 &=(1-c)/2 .
\end{align*}
For a general simplex weight vector $(t,1-t)^\top$, the squared norm becomes
\begin{align*}
 \norm{t g_1+(1-t)g_2}^2
 &= t^2+(1-t)^2-2ct(1-t)\\
 &= 2(1+c)t^2-2(1+c)t+1\\
 &= 2(1+c)(t-1/2)^2+(1-c)/2 .
\end{align*}
The square vanishes at $t=1/2$. Hence as $\theta\downarrow0$ and $c\uparrow1$,
\[
 \chi=\norm d=\sqrt{(1-c)/2}\longrightarrow0,\qquad
 \alpha=1/\sqrt{1-c}\longrightarrow\infty .
\]
The canonical direction still has $\norm v^2=B=2$. Approaching stationarity can therefore make the canonical weights arbitrarily large while the simplex-normalized direction vanishes.

\section{Damped Newton steps for the common weight problem}\label{app:newton}
For a comparison with cyclic coordinate descent, suppose $\chi>0$ and choose $a>0$. Fix a backtracking factor $\omega\in(0,1)$, a fraction $\tau_0\in(0,1)$ of the distance to the positivity boundary, and an Armijo constant $\sigma_N\in(0,1/2)$.

\paragraph{Computing a step.}
Stop as soon as the residual meets the chosen tolerance. Otherwise compute the Newton direction $p$ from
\begin{align*}
 h &= Qa-b/a,\\
 H &= Q+\diag(b_i/a_i^2),\\
 Hp &= -h.
\end{align*}
An exact zero gradient already satisfies the budget equations. For $h\ne0$, start with the largest trial step specified by
\[
 \tau_{\rm pos}=\begin{cases}
 1,& p_i\ge0\text{ for all }i,\\[2pt]
 \min\{1,\tau_0\min_{p_i<0}(-a_i/p_i)\},&\text{otherwise}.
 \end{cases}
\]
Try $\tau=\tau_{\rm pos},\omega\tau_{\rm pos},\omega^2\tau_{\rm pos},\ldots$ in order. Accept the first trial for which
\[
 \Phi_{Q,b}(a+\tau p)\le\Phi_{Q,b}(a)+\sigma_N\tau h^\top p .
\]
Replace $a$ by $a+\tau p$ and repeat.

\paragraph{Why the line search terminates.}
The Hessian is positive definite by Theorem~\ref{thm:exist}. The Newton direction therefore gives
\begin{align*}
 h^\top p &=-h^\top H^{-1}h<0.
\end{align*}
Every trial stays inside the positive orthant. For $p_i\ge0$ we have $a_i+\tau p_i\ge a_i>0$. For $p_i<0$, the definition of $\tau_{\rm pos}$ gives
\begin{align*}
 \tau &\le \tau_0(-a_i/p_i),\\
 \tau p_i &\ge-\tau_0 a_i,\\
 a_i+\tau p_i &\ge (1-\tau_0)a_i>0 .
\end{align*}
By differentiability, the excess over the Armijo threshold satisfies
\[
 \lim_{\tau\downarrow0}
 \frac{\Phi_{Q,b}(a+\tau p)-\Phi_{Q,b}(a)-\sigma_N\tau h^\top p}{\tau}
 =(1-\sigma_N)h^\top p<0 .
\]
Thus every sufficiently small trial passes the test. The trials tend to zero geometrically, so one is accepted after finitely many reductions. Each accepted step decreases the objective and stays in the compact initial sublevel set.

\section{Implementation checklist}\label{app:checklist}
\begin{enumerate}[leftmargin=*]
\item Form $Q$ from the gradients used in the update. When two sets of gradients are used, quantify their discrepancy as in Section~\ref{sec:perturb}.
\item Choose strictly positive budgets $b_i$. They enter the barrier, the coordinate update, and the denominator of the relative residual.
\item Update the positive canonical weights $a$ during a coordinate sweep. Normalize once, when forming the outer mixture $w=a/(\one^\top a)$.
\item Recompute $Qa$ before checking the final residual to remove accumulated error in the cached matrix-vector product.
\item Choose $\delta>0$ for the finite-sweep guarantee of Corollary~\ref{cor:stop}.
\item If an iteration cap is reached while the residual still exceeds $\delta$, report the attained residual and mark the solve as uncertified. At an exactly Pareto stationary point, the positive budget equation has no solution.
\item Report whether direction errors use canonical directions, unit vectors, or simplex-normalized directions. The scale factors can change at every outer iteration.
\item State the stationarity convention alongside the iteration count. The bound $\min_{t<T}\chi(x_t)^2\le C/T$ requires $T\ge C/\varepsilon$ for squared stationarity $\varepsilon$ and $T\ge C/\varepsilon^2$ for stationarity $\varepsilon$.
\end{enumerate}

\end{document}